\documentclass[10pt]{amsart}

\usepackage{latexsym, amsmath,amssymb}

\renewcommand{\epsilon}{\varepsilon}
\newcommand{\newsection}[1]
{\subsection{#1} \setcounter{theorem}{0} \setcounter{equation}{0}
\par\noindent}

\newtheorem{theorem}{Theorem}

\newtheorem{lemma}[theorem]{Lemma}
\newtheorem{corr}[theorem]{Corollary}

\newtheorem{proposition}[theorem]{Proposition}
\newtheorem{deff}[theorem]{Definition}

\newcommand{\bth}{\begin{theorem}}
\newcommand{\ble}{\begin{lemma}}
\newcommand{\bcor}{\begin{corr}}

\newcommand{\bdeff}{\begin{deff}}

\newcommand{\bprop}{\begin{proposition}}
\newcommand{\ele}{\end{lemma}}
\newcommand{\ecor}{\end{corr}}
\newcommand{\edeff}{\end{deff}}

\newcommand{\eprop}{\end{proposition}}

\newcommand{\cd}{\, \cdot\, }

\newcommand{\Rn}{{\mathbb R}^n}

\renewcommand{\Pi}{\varPi}

\renewcommand{\epsilon}{\varepsilon}

\newcommand{\parital}{\partial}
\newcommand{\tidle}{\tilde}

\newcommand{\R}{{\mathbb R}}

\begin{document}

\title[Strichartz-type Estimates for Wave Equation in odd dimension]
{Strichartz-type Estimates for Wave Equation for normally hyperbolic trapped domains}

\thanks{The author would like to thank Prof. Christopher D. Sogge for helpful conversations and suggestions regarding this paper. She would also like to thank Lishui Cheng for his hospitality and encouragement during a recent visit where part of this research carried out.}

\author[H. Sun]{Hongtan Sun}
\address{Department of Mathematics, Johns Hopkins University}

\begin{abstract}
We establish a mixed-norm Strichartz type estimate for the wave equation on Riemannian manifolds $(\Omega,g)$,
for the case that $\Omega$ is the exterior of a smooth, normally hyperbolic trapped obstacle in $n$ dimensional Euclidean space, and $n$ is a positive odd integer.
As for the normally hyperbolic trapped obstacles, we will some loss of derivatives for data in the local energy decay estimate. Hence the global Strichartz estimate has a derivative loss.
However, we can show that the forcing term is bounded by the sum of no more than two Lebesgue $(p,q)$ mixed norms.
\end{abstract}

\maketitle

\newsection{Introduction}

Let $(\Omega,g)$ be a Riemannnian manifold of dimension $n\geq 3$. Throughout this paper, we will assume $\Omega$ to be the region outside a normally hyperbolic trapped obstacle.

The Strichartz-type estimates are a family of space time estimates on the solution 
$u(t,x):(0,T)\times \Omega \longrightarrow \mathbb{C}$
to the wave equation
\begin{equation}\label{1}
\begin{cases}
(\partial^2_t-\Delta_g)u(t,x)=F(t,x), \quad (t,x)\in \R_+\times \Omega
\\
u(0,\cd)=f ,  \quad x\in \Omega 
\\
\partial_t u(0,\cd)=g , \quad x\in \Omega 
\\
u(t,x)=0,\quad \text{on } \, \R_+\times \partial \Omega,
\end{cases}
\end{equation}
where $\Delta_g$ denotes the Laplace-Beltrami operator on $(\Omega,g)$.


If $\Omega$ is a non-trapping domain,
K. Hidano, J. Metcalfe, H. Smith, C. Sogge and Y. Zhou proved in \cite{HMSSZ} the global inhomogeneous mixed-norm Strichartz estimates. The result states that

\begin{equation}\label{2}
\|u\|_{L^p_tL^q_x(\mathbb{R}_{+}\times\Omega)}
\lesssim
\|f\|_{\dot{H}^{\gamma}(\Omega)}
+\|g\|_{\dot{H}^{\gamma}-1(\Omega)}
+\|F\|_{L_t^{r'}L^{s'}_x(\mathbb{R}_{+}\times\Omega)},
\end{equation}
where $\dot{H}^{\gamma}$ denotes the $L^2$ homogeneous Sobolev space over $\Omega$ of order $\gamma$,
and $p,r>2$, $q,s\geq 2$ and they satisfy the gap condition
$$\frac{1}{p}+\frac{n}{q}=\frac{n}{2}-\gamma=\frac{1}{r'}+\frac{n}{s'}-2$$
and the admissible condition
$$\frac{2}{p}+\frac{n-1}{q}\,,\,\frac{2}{r}+\frac{n-1}{s}\leq\frac{n-1}{2}.$$

One of the key inequality to get above \eqref{2} is the plain local energy estimates,
which says
\begin{multline}\label{1.3}
\|u(t,x)\|^2_{L^2_tH^1_x(\mathbb{R}_t\times U)}
+
\|\partial_tu(t,x)\|^2_{L^2_tL^2_x(\mathbb{R}_t\times U)}\\
\lesssim
\|f\|^2_{H^1(V)}
+
\|f\|^2_{L^2(V)}
+\|F(t,x)\|^2_{L^2_tL^2_x(V)},
\end{multline}
for some neighbourhood $U$ and $V$ near the obstacle. 
However, if some of the bicharacteristics on $\Omega$ form closed paths, 
the issue is more intricate. 
In this case, geometrically we say the non-trapping condition fails, and we call such domain a trapped domain.
As a consequence, the plain local energy \eqref{1.3} estimates is not necessarily true \cite{D1} \cite{Bu}.
So the above Strichartz estimates \eqref{2} do not hold for some of the trapped domains,
in particular, the hyperbolic trapped domains.

In fact, the study of wave equations on trapped spaces is fascinating and has attracted many mathematicians. For example, M. Ikawa \cite{IK}, X. Yu \cite{XY}, and H. Christianson \cite{Ch} studied such Strichartz estimate problems for different trapped spaces. And the problem for hyperbolic trapped domains are still open.

In this work, we will establish a Strichartz-type mixed-norm estimates on a hyperbolic tapped domain, 
assuming that the solution satisfies Dirichlet homogeneous boundary conditions. 
The corresponding estimates for problems with Neumann boundary conditions can be established similarly.






\medskip

Let us give a precise description of {\em hyperbolic trapped domain} first. 

Consider $\Omega = \Omega_0\bigsqcup(\mathbb{R}^n\setminus B(0,R))$, 
where $ \Omega_0$ is a smooth compact Riemannian manifold with boundary, 
and the smooth, time-independent metric $g$ equals to the Euclidean metric onb the infinite end $\mathbb{R}^n\setminus B(0,R)$.
So $g_{ij}=\delta_{ij}$ when $|x|>R$. 
For simplicity we assume the compact part $\Omega_0$ and the Euclidean end 
$\mathbb{R}^n\setminus B(0,R)$ are glued together smoothly without any gaps or holes. 
Moreover, we want the trapped set for the geodesic flow on $S^{*}\Omega$ is normally hyperbolic in the following sense, as it is defined in Sec 1.2 of \cite{WZ}.

{\bf Dynamical Assumptions for Normally Hyperbolic Trapped Sets:}

{\em Let $\varphi^t$ denote the geodesic flow, 
and let $r$ denote the distance function to a fixed point in $\Omega$.
Locally define the backward/forward trapped sets by:}

\begin{equation*}
\Gamma_{\pm} = \{\rho\in \pi^{-1} (U_2): \lim_{t\rightarrow\mp\infty} r(\varphi^t(\rho))\neq \infty\}
\end{equation*}

{\em Where $\pi$ is the projection from the cotangent bundles to $\Omega$, and $U_2$ is naborhood near the trapped set $K$, compactly supported in $\Omega$.}

{\em Then the trapped set is defined to be }

$$K=\Gamma_{+}\cap\Gamma{-}$$

{\em Let $p$ be the principle symbol of $\Delta_g$, and write $K_{\lambda}=K\cap p^{-1}(\lambda)$.
Then the} {\bf Dynamical Hypotheses} {\em are }

{\em 1, There exists $\delta>0$ such that $dp\neq 0$ on $p^{-1}(\lambda)$ for $|\lambda|<\delta$.}

{\em 2, $\Gamma{\pm}$ are codimension-one smooth manifolds intersecting transversely at $K$.}

{\em 3, The flow is hyperbolic in the normal direction to $K$ within the energy surface: 
there exists subbundles $E^{\pm}$ of $T_{K_{\lambda}}\Gamma^{\lambda}_{\pm}$ such that}

$$T_{K_{\lambda}}\Gamma^{\lambda}_{\pm}=TK_{\lambda}\oplus E^{\pm}$$

{\em where }

$$d\varphi^t : E^{\pm}\rightarrow E^{\pm}$$

{\em and there exists $\theta>0$ such that for all $|\lambda|<\delta$,}

$$d\varphi^t(v)\leq Ce^{-\theta |t|}||v||$$

{\em for all $v\in E^{\mp}$, $\pm t \geq 0$.}

\medskip

Now we can consider the Laplace-Beltrami operator  $\Delta_g$  associated with $g$,
We shall assume throughout that the spatial dimension $n\geq 3$ is odd.


In order to describe the Sobolev spaces that is used on $\Omega$, 
we shall start with the regular $L^2$ Sobolev space on $\mathbb{R}^n$. 
The definition is the same as the definition of the Sobolev spaces used in \cite{HMSSZ} and \cite{XY}.

Recall that $\dot{H}^{\gamma}(\mathbb{R}^n)$ is the $L^2$ homogeneous Sobolev space with norm given by

\begin{equation*}
\|f\|^2_{\dot{H}^{\gamma}(\mathbb{R}^n)}
=
\|(\sqrt{-\triangle})^{\gamma}f\|^2_{L^2(\mathbb{R}^n)}
=
(2\pi)^{-n}\int_{\mathbb{R}^n}||\xi|^{\gamma}\hat{f}(\xi)|^2 d\xi.
\end{equation*}

While the inhomogeneous Sobolev space $H^{\gamma}(\mathbb{R}^n)$ has norm defined by

\begin{equation*}
\|f\|^2_{H^{\gamma}(\mathbb{R}^n)}
=
\|(1-\triangle)^{\gamma /2}f\|^2_{L^2(\mathbb{R}^n)}
=
(2\pi)^{-n}\int_{\mathbb{R}^n}|(1+|\xi|^2)^{\gamma /2}\hat{f}(\xi)|^2 d\xi,
\end{equation*}

with $\hat{f}$ denoting the Fourier transform and $\triangle$ denoting the standard Laplacian.

Then take a smooth cut-off function $\beta$ on $\R^n$. 
Assume that $\beta$ is supported in $|x|<2R$, 
and $\beta(x)=1$ when $|x|<R$.
Let $\Omega'$ be the embedding
of $\Omega\cap\{|x|<2R\}$ into  a compact manifold with boundary,
so that $\partial\Omega'=\partial\Omega$. 
Then we can define the inhomogeneous $L^2$ Sobolev norm on $\Omega$

\begin{equation}\label{1.8}
\|f\|_{H^{\gamma}(\Omega)}=\|\beta f\|_{H^{\gamma}(\Omega')}
+\|(1-\beta)f\|_{H^{\gamma}(\mathbb{R}^n)},
\end{equation}

and the homogeneous $L^2$ Sobolev norm on $\Omega$

\begin{equation}\label{1.10}
\|f\|_{\dot{H}^{\gamma}(\Omega)}=\|\beta f\|_{H^{\gamma}(\Omega')}
+\|(1-\beta)f\|_{\dot{H}^{\gamma}(\mathbb{R}^n)}.
\end{equation}

Notice that the spaces $H^{\gamma}(\Omega')$ are defined by a spectral decomposition of $\Delta_g|_{\Omega'}$, subject to the Diriclet boundary condition.
In fact, $H^{-\gamma}(\Omega)$ is the dual of $H^{\gamma}(\Omega)$, 
and $\dot{H}^{-\gamma}(\Omega)$ is dual to $\dot{H}^{\gamma}(\Omega)$.

As we mentioned at the beginning, 
the plain local energy estimates fails in hyperbolic trapped case.
However, later we will see that we can get a similar estimates with some derivative loss.
And this new local energy estimates helped us to conquer some of the difficulties.
And in order to deal with the extra derivative loss in the local energy estimate, 
we introduce a Sobolev-type norm in the follow, as it did in \cite{XY}.

\begin{deff}
Define $\tilde{H}^{\gamma}_{\varepsilon}(\mathbb{R}^n)$ 
 to be the space with norm

\begin{equation}\label{1.5}
\|h\|_{\tilde{H}^{\gamma}_{\varepsilon}(\mathbb{R}^n)}
=
\||D|^{\gamma}(1-\Delta)^{\varepsilon /2}h\|_{L^2_x(\mathbb{R}^n)}
=
(\int_{\mathbb{R}^n}||\xi|^{\gamma}(1+|\xi|^2)^{\varepsilon /2}\hat{h}(\xi)|^2 d\xi)^{1/2}
\end{equation}

\end{deff}

Notice that, for both positive and negative $\varepsilon$, 
above Sobolev-type norm has the following properties:

\begin{equation}\label{1.6}
\|h\|_{\tilde{H}^{\gamma}_{\varepsilon}(\mathbb{R}^n)}
\approx
\|h\|_{\dot{H}^{\gamma}(|\xi|<1)}
+
\|h\|_{\dot{H}^{\gamma+\varepsilon}(|\xi|>1)}.
\end{equation}

And for non-negative $\varepsilon$, we even have
\begin{equation}\label{1.9}
\|h\|_{\tilde{H}^{\gamma}_{\varepsilon}(\mathbb{R}^n)}
\approx
\|h\|_{\dot{H}^{\gamma}(\Rn)}
+
\|h\|_{\dot{H}^{\gamma+\varepsilon}(\Rn)}.
\end{equation}

Now we can define the Sobolev-type norm on a manifold $\Omega$ similarly, 
as what we did earlier in \eqref{1.8} and \eqref{1.10}.

\begin{equation}\label{1.7}
\|f\|_{\tilde{H}^{\gamma}_{\varepsilon}(\Omega)}
=
\|\beta f\|_{\tilde{H}^{\gamma}_{\varepsilon}(\Omega')}
+
\|(1-\beta) f\|_{\tilde{H}^{\gamma}_{\varepsilon}(\mathbb{R}^n)}
\end{equation}

When $\varepsilon=0$, our definition turns out to coincide with the homogeneous Sobolev spaces, 
and we have

\begin{equation*}
{\tilde{H}^{\gamma}_0(\mathbb{R}^n)}=\dot{H}^{\gamma}(\mathbb{R}^n),
\quad
{\tilde{H}^{\gamma}_0(\Omega)}=\dot{H}^{\gamma}(\Omega)
\end{equation*}

\medskip

Now let us consider the wave equation \eqref{1} on hyperbolic trapped domain $(\Omega,g)$.
We assume $n$ is odd, $\varepsilon>0$,
$p>2$, $\gamma\in (-\frac{n-3}{2},\frac{n-1}{2})$,and $(p,q,\gamma)$ satisfies

\begin{equation}\label{1.2}
\frac{1}{p}+\frac{n}{q}=\frac{n}{2}-\gamma,\quad
\begin{cases}
\frac{3}{p}+\frac{2}{q}\leq 1, \quad n=3
\\
\frac{2}{p}+\frac{2}{q}\leq 1, \quad n>3
\end{cases}
\end{equation}

Our first result concerning Strichartz-type mixed-norm estimates is the following.

\begin{theorem}\label{thm1}
Let $n\geq 3$ be odd, 
and fix $\Omega=\Omega_0\cup(\mathbb{R}^n\setminus B(0,R))$ 
as a normally hyperbolic trapped domain, and $(p,q,\gamma)$
satisfying above \eqref{1.2}.
Assume that $p>2$ and $\gamma\in(-\frac{n-3}{2},\frac{n-1}{2})$.

Let $F=0$ and $u=u(t,x)$ solves \eqref{1},
then for all $\varepsilon>0$,

\begin{equation}\label{1.1}
\|u\|_{L^p_tL^q_x(\mathbb{R}_{+}\times \Omega)}
\lesssim
\|f\|_{\tilde{H}^\gamma_{\varepsilon}(\Omega)}
+
\|g\|_{\tilde{H}^{\gamma-1}_{\varepsilon}(\Omega)}.
\end{equation}
\end{theorem}

Theorem \ref{thm1} gives the global mixed-norm Stricharts-type estimates for solutions of homogeneous wave equation on $\Omega$.
As we mentioned, there are derivative losses in the local energy estimates. 
And \eqref{1.1} shows that the global mixed-norm estimates inherited those derivative losses.
However, later in the proof we can see 
that the homogeneous local Strichartz-type mixed-norm estimate does not have such derivative loss.

The appearance of derivative loss in the global estimates 
makes it hard to get a good estimates for the inhomogeneous equation.
For example, X. Yu encountered similar difficulties in \cite{XY}. 
So she has to construct a certain weighted-norm Strichartz-type estimates.

In this work, we get through this difficulty 
by showing that we need no more than two Lebesgue $(p,q)$ mixed-norms
to bound the inhomogeneous solution. 
So the main Strichartz-type estimate is the following.

\begin{theorem}\label{thm2}
Let $n\geq 3$ be odd, and fix a normally hyperbolic trapped domain as it is in Theorem \ref{thm1} above. Assume $p>2$,  $\gamma \in (-\frac{n-3}{2},\frac{n-1}{2})$, and $(p,q,\gamma)$ and $(r, s, 1-\gamma)$ both satisfy \eqref{1.2}.

Let $u=u(t,x)$ solve \eqref{1}.
 Then for any $\varepsilon, \delta>0$, we have

\begin{multline}\label{1.4}
\|u\|_{L^p_tL^q_x(\mathbb{R}_{+}\times \Omega)}
\lesssim
\|f\|_{\tilde{H}^\gamma_{\varepsilon}(\Omega)}
+
\|g\|_{\tilde{H}^{\gamma-1}_{\varepsilon}(\Omega)}
\\
+
\|F\|_{L^{r'}_tL^{s'}_x(\mathbb{R}_{+}\times\Omega)}
+
\|F\|_{L^{r'}_tL^{s'-\delta}_x(\mathbb{R}_{+}\times\Omega)}.
\end{multline}
\end{theorem}

\medskip

The idea and method of the proofs are inherited from \cite{HMSSZ} and \cite{XY}. 
In section $2$, we will prove the afore mentioned local energy estimate. It is different from the well-know plain local energy estimates since it contains some derivative loss, which comes out naturally from the geometry of the normally hyperbolic trapping. In section $3$, we show a local Strichartz-type mixed-norm estimate for homogeneous wave equation on $\Omega$. Even considering the existence of trapped geodesics, the local estimates still looks the same as the usual one. However, careful readers may notice that the domain for the Lebesgue index $(p,q)$ here is much smaller than those for the Minkowski case. In section $4$, we combine the well-known global Strichartz estimate for wave equations on Minkowski space with the results in section $2$ and $3$, and get the global Strichartz type estimate and prove theorem \ref{thm1}. This part of proof is a direct application of a result of X. Yu \cite{XY}.  In section $5$, we prove the inhomogeneous estimates in theorem \ref{thm2} using a $TT^{*}$ argument.

\qed


\newsection{Local Energy Decay Effect for Weakly Trapped Domain}

Consider the normally hyperbolic trapped domain$\Omega = \Omega_0\bigsqcup(\mathbb{R}^n\setminus B(0,R))$ described in Section 1, satisfy the afore mentioned Dynamical Assumptions. Assume that $u=u(t,x)$ solves \eqref{1} with $f$, $g$, $F$ compactly supported near $\partial\Omega$.

Then a special case of Theorem 3 in Section 5 of \cite{WZ} is,

\begin{corr}
Supposed that $\Omega$ is a normally hyperbolic trapped domain as above.
Assume that, if for a fixed $R>0$, $u$ solves \eqref{1}. Assume that $F=0$ and $f$,$g$ compactly supported in $\{|x|<R\}$. Then there exists $\alpha>0$, $K\in \mathbb{Z}^{+}$, 
for some $C=C(\alpha, K, R)$, we have
\begin{equation}\label{2.4}
\int_{|x|<R}(|u'(t,x)|^2+|\partial_tu(t,x)|^2)dx
\leq
Ce^{-\alpha t}(||f||^2_{H^{K+1}(|x|<R)}+||g||^2_{H^K(|x|<R)})
\end{equation}
\end{corr}

Interpolate \eqref{2.4} with the energy estimate, we get the following lemma.

\begin{lemma}
Fix $R>0$. If $u=u(t,x)$ solves \eqref{1} with $f$, $g$, $F$ compactly supported near $\partial\Omega$, $\exists c>0$, so that
\begin{equation}\label{2.1}
\|u'(t,x)\|_{L^2_x(|x|<R)} \lesssim e^{-ct}
(
\|f\|_{H^{1+\varepsilon}(|x|<R)}
+
\|g\|_{\dot{H}^{\varepsilon}(|x|<R)}
)
\end{equation}
For any $\varepsilon>0$
\end{lemma}

Lemma 2.2 is a local energy decay estimate. It states that the local energy near the trapped sets decays exponentially with respect to the time variable $t$, comparing to the Soblove norms of the initial data. Using this estimate, we can get the local energy estimate, as it is stated in the following proposition.

\begin{proposition}\label{prop2.3}
Assume $u=u(t,x)$ solves \eqref{1}, and $f$, $g$, $F$ vanishes for $|x|>R$. Then we have
\begin{multline}\label{2.2}
\int_0^{\infty}(\|u(t,x)\|^2_{\dot{H}^1(|x|<R)}+\|\partial_tu(t,x)\|^2_{L^2(|x|<R)}) dt
\\
\lesssim \|f\|^2_{\dot{H}^{1+\varepsilon}(\Omega)}
+\|g\|^2_{\dot{H}^{\varepsilon}(\Omega)}
+\int_0^{\infty}\|F(s,x)\|^2_{\dot{H}^{\varepsilon}(\Omega)}ds.
\end{multline}
\end{proposition}

Now let us prove Proposition \ref{prop2.3}. Notice that the functions of both sides of the equation are compactly supported, the left hand side of \eqref{2.2} is equivalent to

\begin{equation}\label{a2}
A=\|u'(t,x)\|^2_{L^{2}_tL^2_x(\mathbb{R}_{+}\times\{|x|<R\})}.
\end{equation}

When $F=0$, taking square and integral on both sides of \eqref{2.1},
we get the desired homogeneous estimates,

\begin{equation}\label{A2}
A^2 \lesssim \|f\|^2_{\dot{H}^{1+\varepsilon}(\Omega)}
+
\|g\|^2_{\dot{H}^{\varepsilon}(\Omega)}.
\end{equation}


Assume that $D_g=\sqrt{-\Delta_g}$. Taking $F=0$, we can get from \eqref{2.1} that 

\begin{equation}\label{2.3}
\|e^{it|D_g|}h \|_{L^2_x(|x|<R)}
\lesssim
e^{-ct}\|h\|_{\dot{H}^{\varepsilon}(\Omega)},
\end{equation}

for $h \in \dot{H}^{\varepsilon}(\Omega)$ Supported in ${|x|<R}$.

Then by Duhammel's Principle and applying Minkowski inequality to 
$$\|e^{it|D_g|}h \|_{L^2_x(|x|<R)}
\lesssim
e^{-ct}\|h\|_{\dot{H}^{\varepsilon /2}(\Omega)},$$
we get the desired estimate for $F\neq 0$.

\qed

\medskip

\newsection{Local Strichartz Estimates}

In order to prove the local Strichartz estimates,
we will need a special case of a lemma of Christ and Kiselev \cite{CK}.
So we quote it here as the following lemma. 
A proof of Christ-Kiselev Lemma could be found in \cite{So}.

\begin{lemma}$(${\bf Christ-Kiselev Lemma}$)$
Let $X$ and $Y$ be Banach Spaces and assume that $K(t,s)$ is a continuous function taking its values in $B(X,Y)$, the space of bounded linear mappings from $X$ to $Y$. Suppose that 
$-\infty\leq a < b \leq \infty$, and set
\begin{equation*}
Tf(t)=\int_a^bK(t,s)f(s) ds 
\end{equation*}.
Assume that 
\begin{equation}
\|Tf\|_{L^q([a,b];Y)}\leq \|f\|_{L^p([a,b];Y)}
\end{equation}.
Set
\begin{equation*}
Wf(t)=\int_t^bK(t,s)f(s) ds 
\end{equation*}.
Then if $1\leq p<q\leq \infty$,
\begin{equation}
\|Wf\|_{L^q([a,b];Y)}
\leq
\frac{2^{-2(1/p-1/q)}\cdot 2C}{1-2^{-(1/p-1/q)}}
\|f\|_{L^p([a,b];X)}
\end{equation}.
\end{lemma}

\medskip

The main result of this section is the following proposition:

\begin{proposition}\label{prop3.2}
If $u=u(t,x)$ solves \eqref{1} with $F=0$. Assume that $p>2$
and $\gamma \in [-\frac{n-3}{2},\frac{n-1}{2})$. Then for $(p,q,\gamma)$ satisfying \eqref{1.2}, we have, for all $\varepsilon>0$,

\begin{equation}\label{3.1}
\|u\|_{L^p_tL^q_x([0,1]\times \Omega)}
\lesssim
\|f(x)\|_{\dot{H}_{\gamma}(\Omega)}
+
\|g(x)\|_{\dot{H}^{\gamma-1}(\Omega)}
\end{equation}

\end{proposition}

In order to prove Proposition \ref{prop3.2}, let us consider a smooth cutoff function $\varphi\in C_0^\infty(\mathbb{R}^n)$. 
Assume that $\varphi(x)=0$ when $|x|>2R$, and $\varphi(x)=1$ when $|x|<R$.

Then for the solution $u$ in Proposition \ref{prop3.2}, by taking $v=\varphi u$ and $w=(1-\varphi)u$, 
we can write it as

\begin{equation}
u=v+w=\varphi u+(1-\varphi)u
\end{equation}

\medskip

The second part, $w=(1-\varphi)u$, solves the following equation, 
as $g_{ij}=\delta_{ij}$ when $|x|>R$.

\begin{equation}\label{3.2}
\begin{cases}
(\partial^2_t-\Delta)w=-[\varphi,\Delta_g]u
\\
w(0,x)=(1-\varphi(x))f(x)
\\
\partial_tw(0,x)=(1-\varphi(x))g(x)
\end{cases}
\end{equation}

Notice that, as $w$ vanishes in $|x|<R$, 
and $G(t,x)=-[\varphi,\Delta_g]u$ is supported on $\{R<|x|<2R\}$, 
\eqref{3.2} is a wave equation in the Minkowski Space.
So if $w_0$ solves the corresponding homogeneous equation 
$(\partial^2_t-\Delta)w_0=0$ ,
by Corollary 1.2 in chapter 4 of \cite{So},
we have,

\begin{equation}\label{3.6}
\|w_0\|_{L^p_tL^q_x([0,1]\times \mathbb{R}^n)}
\lesssim
\|(1-\varphi(x))f(x)\|_{\dot{H}^{\gamma}(\mathbb{R}^n)}
+
\|(1-\varphi(x))g(x)\|_{\dot{H}^{\gamma-1}(\mathbb{R}^n)}
\end{equation}

for $(p,q, \gamma)$ satisfying \eqref{1.2}.

On the other hand, if $w_1$ solves the inhomogeneous equation
$(\parital_t-\Delta)w_1=F(t,x)$  on Minkowski space, with zero initial data. 
By Duhammel's Principle,

\begin{multline}\label{3.3}
\|w_1\|_{L^p_tL^q_x([0,1]\times \Rn)}
\lesssim
\int_0^1\|\frac{sin(t-s)\sqrt{-\Delta}}{\sqrt{-\Delta}}F(s,\cdot)\|_{L^p_tL^q_x([0,1]\times \Rn)} ds
\\
\lesssim \|F(t,x)\|_{L^1_t\dot{H}^{\gamma-1}([0,1]\times\Rn)}
\lesssim \|F(t,x)\|_{L^2_t\dot{H}^{\gamma-1}([0,1]\times\Rn)}
\end{multline}

So combine \eqref{3.6} and \eqref{3.3}, we get, 

\begin{multline}
\|w\|_{L_t^pL_x^q([0,1]\times \Omega)}
\lesssim
\|(1-\varphi)f\|_{\dot{H}^\gamma(\mathbb{R}^n)}
+\|(1-\varphi)g\|_{\dot{H}^{\gamma-1}(\mathbb{R}^n)}
\\
+\|[\varphi,\Delta_g]u\|_{L^2_t\dot{H}^{\gamma-1}([0,1]\times\mathbb{R}^n)}
\end{multline}

As 
$[\varphi,\Delta_g]u
=
-(\Delta_g\varphi)u
+
\nabla_g\varphi\cdot\nabla_gu$, supported on $R<|x|<2R$, according to the local energy effect, we can get

\begin{equation}
\|[\varphi,\Delta_g]u\|_{L^2_t\dot{H}^{\gamma-1}([0,1]\times\mathbb{R}^n)}
\lesssim
\|\varphi f\|_{\dot{H}^{\gamma}(\Omega)}
+
\|\varphi g\|_{\dot{H}^{\gamma-1}(\Omega)}
\end{equation}

So $w$ is bounded by
\begin{equation*}
\|(1-\varphi)f\|_{\dot{H}^\gamma(\mathbb{R}^n)}+
\|(1-\varphi)g\|_{\dot{H}^{\gamma-1}(\mathbb{R}^n)}+
\|\varphi f\|_{\dot{H}^{\gamma}(\Omega)}+
\|\varphi g\|_{\dot{H}^{\gamma+-1}(\Omega)}
\end{equation*}

Which is equivalent to 
$$\|f\|_{\dot{H}^{\gamma}(\Omega)}
+
\|g\|_{\dot{H}^{\gamma-1}(\Omega)}$$

\medskip

Now let us consider $v=\varphi u$.

Notice that $v$ satisfies the following wave equation,

\begin{equation}\label{3.5}
\begin{cases}
(\partial^2_t-\Delta_g)v=[\varphi,\Delta_g]u
\\
v(0,x)=\varphi(x)f(x)
\\
\partial_tv(0,x)=\varphi(x)g(x)
\\
v(t,x)|_{\partial\Omega}=0
\end{cases}
\end{equation}

Since $v$ itself, the forcing term and the initial condition 
are all compactly supported in $\{x \in \Omega:|x|<2R\}$,
 we can consider the wave equation as on a compact manifold $(\tilde{\Omega},\tilde{g})$ with boundary.
 We can take $\tilde{\Omega}$ so that $\partial\tilde{\Omega}=\partial\Omega$, and $\tilde{g}$ coincide with $g$.

Then by Theorem 1.1 in \cite{BSS}, if $(p,q,\gamma)$ satisfies \eqref{1.2}, and $v_0$ solves the corresponding homogeneous wave equation,
\begin{equation}\label{3.7}
\begin{cases}
(\partial^2_t-\Delta_g)v_0=0
\\
v_0(0,x)=\varphi(x)f(x)
\\
\partial_tv_0(0,x)=\varphi(x)g(x)
\\
v_0(t,x)|_{\partial\Omega}=0
\end{cases}
\end{equation}

we have

\begin{equation}\label{3.10}
\|v_0\|_{L_t^pL_x^q([0,1]\times \Omega)}\lesssim
\|\varphi f\|_{\dot{H}^\gamma(\tilde{\Omega})}
+\|\varphi g\|_{\dot{H}^{\gamma-1}(\tilde{\Omega})}
\end{equation}

By Duhammel's principle and the spectrum decomposition of the Sobolev norm on compact manifolds, we get that 

\begin{equation}\label{3.8}
\|v_1\|_{L_t^pL_x^q([0,1]\times \Omega)}
\lesssim
\|F(t,x)\|_{L^2_tH^{\gamma-1}([0,1]\times\tilde{\Omega})}
\end{equation}

If $v_1$ solves
\begin{equation}
\begin{cases}
(\partial^2_t-\Delta_g)v_1=F(t,x) \quad \text{on} \quad [0,1]\times\tilde{\Omega}
\\
v_1(0,x)=0
\\
\partial_tv_1(0,x)=0
\\
v_1(t,x)|_{\partial\Omega}=0
\end{cases}
\end{equation}

So by \eqref{3.10} and \eqref{3.8}, it follows that
\begin{equation}
\|v\|_{L_t^pL_x^q([0,1]\times \Omega)}
\lesssim
\|\varphi f\|_{\dot{H}^\gamma(\tilde{\Omega})}
+\|\varphi g\|_{\dot{H}^{\gamma-1}(\tilde{\Omega})}
+\|[\varphi,\Delta_g]u\|_{L^2_tH^{\gamma-1}([0,1]\times\tilde{\Omega})}
\end{equation}

where the right hand side of the inequality is bounded by 
$\|f\|_{\dot{H}^{\gamma}(\Omega)}+
\|g\|_{\dot{H}^{\gamma-1}(\Omega)}$, as
$\|[\varphi,\Delta_g]u\|_{L^2_tH^{\gamma-1}([0,1]\times\tilde{\Omega})}
\approx
\|[\varphi,\Delta_g]u\|_{L^2_tH^{\gamma-1}([0,1]\times\Omega)}$

Therefore, by Minkowski's inequality, 
we got the desired local Strichartz estimates.

\qed

\medskip


\newsection{Global Strichartz Estimates for Homogeneous Wave Equation}

There is a long history to establishing the global Strichartz estimate in Minkowski space, beginning with the original work by Strichartz \cite{Str}. Some subsequent work is done by Genibre-Velo \cite{GV}, Pecher \cite{Pe}, Kapitanski \cite{Ka}, Lindblad-Sogge \cite{LS} \cite{LS2}, Mockenhaupt-Seeger-Sogge \cite{MSS}, Keel-Tao \cite{KT}, etc. 
One of the result is the following estimate, which can be found as Corollary 2.1 in chapter 4 of \cite{So}.

\noindent{\bf Global Minkowski Strichartz estimates.}

For solutions to \eqref{1} with $F=0$, in the case of $\Omega=\mathbb{R}^n$ and $g_{ij}=\delta_{ij}$, 
the following holds, when $p>2$, $(p,q,\gamma)$ satisfies \eqref{1.2}, we have that,

\begin{equation}\label{4.1}
\|u\|_{L^p_tL^q_x(\mathbb{R}^{1+n})}\lesssim
\|f\|_{\dot{H}^\gamma(\mathbb{R}^n)}
+\|g\|_{\dot{H}^{\gamma-1}(\mathbb{R}^n)}
\end{equation}

Xin Yu defined {\bf almost admissibility} in \cite{XY}, we state the definition here,

\begin{deff}
We say that $(X,\gamma,\eta,p)$ is almost admissible if it satisfies

i), Minkowski almost Strichartz estimates
\begin{equation}
\|u\|_{L^p_tX([0,S]\times\mathbb{R}^n)}
\lesssim
\|u(0,\cd)\|_{\dot{H}^\gamma(\mathbb{R}^n)}
+\|\partial_tu(0,\cd)\|_{\dot{H}^{\gamma-1}(\mathbb{R}^n)}
\end{equation}

ii), Local almost Strichartz estimates for $\Omega$
\begin{equation}
\|u\|_{L^p_tX([0,1]\times\Omega)}
\lesssim
\|u(0,\cd)\|_{\tidle{H}_\eta^\gamma(\Omega)}
+\|\partial_tu(0,\cd)\|_{\tidle{H}_\eta^{\gamma-1}(\Omega)}
\end{equation}
\end{deff}

Xin Yu \cite{XY} also proved the following Generalized Strichartz Estimates Theorem.

\begin{theorem}[Theorem 1.5 in \cite{XY} for Dirichlet-wave equation]
Let $n>2$ and assume that $(X,\gamma,\eta,p)$ is almost admissible with
$$p>2\quad and \quad \gamma\in( -\frac{n-3}{2},\frac{n-1}{2}).$$
Then if the local smoothing estimate \eqref{2.2} is valid and 
if u solves \eqref{1} with forcing term $F=0$,
we have the abstract Strichartz estimates
\begin{equation}
\|u\|_{L^p_tX([0,\infty)\times\Omega}
\lesssim
\|f\|_{\tilde{H}^{\gamma}_{\varepsilon+\eta}(\Omega)}
+
\|g\|_{\tilde{H}^{\gamma-1}_{\varepsilon+\eta}(\Omega)}
\end{equation}
\end{theorem}
Now from proposition \ref{prop3.2} and \eqref{4.1}, we can see that $(L^q(\Omega),\gamma, 0, p)$ is almost admissible, if $\gamma\in(-\frac{n-3}{2},\frac{n-1}{2})$, $(p,q,\gamma)$ satisfying \eqref{1.2} and any $\varepsilon>0$. So Theorem \ref{thm1} follows by applying Theorem 4.2.
\qed

\medskip

\newsection{Inhomogeneous Estimates}

Notice that the homogeneous Strichartz estimates in Theorem \ref{thm1} is different from the usual Strichartz estimates we have seen for Euclidean space, non-trapping domain or manifolds with or without boundary.
Because it contains derivative loss on the right hand side of the inequality. 
This cause problems for getting a inhomogeneous estimate that could maintain the order of derivatives on both sides of the inequality. 
So when dealing with this problem, 
I need to use two different Lebesgue $(p,q)$ mixed-norms to hedge the potential derivative loss.

\noindent{\bf Proof of Theorem \ref{thm2}:}

By theorem \ref{thm1}, it suffices to prove the case for $f=g=0$.

By Duhammel's Principle, we have to show,

\begin{multline}\label{5.1}
\|\int_0^te^{i(t-s)|D|}|D|^{-1}F(s,\cdot) ds\|_{L^p_tL^q_x(\mathbb{R}_{+}\times \Omega)}
\\
\lesssim
\|F(t,x)\|_{L^{r'}_tL^{s'}_x(\mathbb{R}_{+}\times \Omega)}
+
\|F(t,x)\|_{L^{r'}_tL^{s'-\delta}_x(\mathbb{R}_{+}\times \Omega)},
\end{multline}

for $(p,q,\gamma)$ and $(r,s,1-\gamma)$ satisfying \eqref{1.2} respectively, and any $\delta>0$

By Christ-Kiselev Lemma, it suffices to show

\begin{multline}\label{5.2}
\|\int_0^{\infty}e^{i(t-s)|D|}|D|^{-1}F(s,\cdot) ds\|_{L^p_tL^q_x(\mathbb{R}_{+}\times \Omega)}
\\
\lesssim
\|F(t,x)\|_{L^{r'}_tL^{s'}_x(\mathbb{R}_{+}\times \Omega)}
+
\|F(t,x)\|_{L^{r'}_tL^{s'-\delta}_x(\mathbb{R}_{+}\times \Omega)}
\end{multline}

From theorem \ref{thm1}, we can see the left hand side of \eqref{5.2} is bounded by

\begin{equation}\label{5.3}
\|\int_0^{\infty}e^{-is|D|}|D|^{-1}F(s,\cdot) ds\|_{\tilde{H}^{\gamma}_{\varepsilon}(\Omega)},
\end{equation}

for any $\varepsilon>0$.

Let $P=\sqrt{-\Delta_g}$. According to \eqref{1.6}, this is equals to

\begin{equation}\label{5.4}
\|\int_0^{\infty}e^{-is|D|}F(s,\cdot) ds\|_{\dot{H}^{\gamma-1}(|P|\leq 1)}
+
\|\int_0^{\infty}e^{-is|D|}F(s,\cdot) ds\|_{\dot{H}^{\gamma+\varepsilon-1}(|P|\geq 1)}
\end{equation}

If $(r,s,1-\gamma)$ satisfies \eqref{1.2},
for $\varepsilon>0$,
the duality of \eqref{1.1} for $(r, s, 1-\gamma, 0)$ gives

\begin{equation}\label{5.5}
\|\int_0^{\infty}e^{-is|D|}F(s,\cdot) ds\|_{\tilde{H}^{\gamma-1}_{-\varepsilon}(\Omega)}
\lesssim
\|F(t,x)\|_{L^{r'}_tL^{s'}_x(\mathbb{R}_{+}\times\Omega)}
\end{equation}

And by definition of the $\tilde{H}^{\gamma}_{\varepsilon}$, 
the left hand side of \eqref{5.5} is equivalent to 

\begin{equation}\label{5.6}
\|\int_0^{\infty}e^{-is|D|}F(s,\cdot) ds\|_{\dot{H}^{\gamma-1}(|P|\leq 1)}
+
\|\int_0^{\infty}e^{-is|D|}F(s,\cdot) ds\|_{\dot{H}^{\gamma-\varepsilon-1}(|P|\geq 1)}
\end{equation}

So the first term in \eqref{5.4} is bounded by 
$\|F(t,x)\|_{L^{r'}_tL^{s'}_x(\mathbb{R}_{+}\times \Omega)}$.
And this is the first term on the left hand side of \eqref{5.1}.

Similarly, consider $(r,\tilde{s},1-(\gamma+2\varepsilon))$ satisfying \eqref{1.2}. 
Then the duality of \eqref{1.1} for $(r, s, 1-(\gamma+2\varepsilon), 0)$ gives

\begin{equation}\label{5.7}
\|\int_0^{\infty}e^{-is|D|}F(s,\cdot) ds\|_{\tilde{H}^{\gamma+2\varepsilon-1}_{-\varepsilon}(\Omega)}
\lesssim
\|F(t,x)\|_{L^{r'}_tL^{\tilde{s}'}_x(\mathbb{R}_{+}\times\Omega)}
\end{equation}

And the left and side of \eqref{5.7} turns out to be,

\begin{multline}\label{5.8}
\|\int_0^{\infty}e^{-is|D|}F(s,\cdot) ds\|_{\tilde{H}^{\gamma+2\varepsilon-1}_{-\varepsilon}(\Omega)}
\approx
\|\int_0^{\infty}e^{-is|D|}F(s,\cdot) ds\|_{\dot{H}^{\gamma+\varepsilon-1}(|P|\geq 1)}
\\
+
\|\int_0^{\infty}e^{-is|D|}F(s,\cdot) ds\|_{\dot{H}^{\gamma+2\varepsilon-1}(|P|\leq 1)}
\end{multline}

Which provides the bound for the second term in \eqref{5.4}. 
Notice that $\varepsilon$ can be arbitrarily small, 
it follows that we can choose $\tilde{s}'=s'-\delta$,
for arbitrarily small $\delta>0$. And this completed the proof for Theorem \ref{thm2}.

\qed


\begin{thebibliography}{MA}



\bibitem{BSS} M. Blair, H. Smith and C. D. Sogge:
{\em Strichartz estimates for the wave equation on manifolds with boundary},
arXiv:0805.4733.

\bibitem{Burq} N. Burq: {\em Global Strichartz estimates for
nontrapping geometries: About an article by H. Smith and C. Sogge},
Comm. Partial Differential Equations  {\bf 28} (2003), 1675--1683.

\bibitem{Bu} N. Burq :{\em Smoothing effect for Schr\"{o}dinger boundary value problems},
Duke Mathematical Journal, Vol. 123, No. 23 (2004), 403-427.

\bibitem{BGH} N. Burq, C. Guillarmou, A. Hassell:
{\em Strichartz Estimates without Loss on Manifolds with Hyperbolic Trapped Geodesics},
Geom. Funct. Anal. Vol. 20(2010) 627-656

\bibitem{CK} M. Christ and A. Kiselev:
{\em Maximal functions associated to filtrations},
J. Funct. Anal. {\bf 179} (2001), 409--425.

\bibitem{Ch} H. Christianson,
{\em Applications of Cutoff Resolvent Estimates to the Wave Equation},
2007, Arxiv 0709.0555

\bibitem{CW} H. Christianson and J. Wunsch,
{\em Local Smoothing for the Schr\"{o}dinger Equation with a Prescribed Loss}
(2012) arXiv:1103.3908

\bibitem{D1}S. I. Doi:
{\em Smoothing Effects of Schr\"{o}dinger Evolution Groups on Riemannian Manifolds},
Duke Math J. {\bf 82}(1996), 679-706. MR 1387689

\bibitem{D2}S. I. Doi:
{\em Smoothing Effects for Schr\"{o}dinger Evolution Equations and Global Behaviour of Geodesic Flow},
Math. Ann. {\bf 318}(2000), 355-389. MR 1795567

\bibitem{DMSZ}Y. Du, J. Metcalfe, C. D. Sogge and Y. Zhou:
{\em Concerning the Strauss conjecture and almost global existence
for nonlinear Dirichlet-wave equations in $4$-dimensions},
Comm. Partial Differential Equations {\bf 33} (2008), 1487--1506.

\bibitem{FW} D. Fang and C. Wang:
{\em Some Remarks on Strichartz Estimates for Homogeneous Wave Equation},
Nonlinear Analysis {\em 65}(2006) 697-706

\bibitem{GLS} V. Georgiev, H. Lindblad, C. D. Sogge:
{\em Weighted Strichartz Estimates and Global Existence for Semilinear Wave Equations},
American Journal of Mathematics {\bf 119}(1997), 1291-1319

\bibitem{GV} J. Ginibre, G. Velo:
{\em Generalized Strichartz Inequalities for the Wave Equation}
Journal of Functional Analysis {\bf 133.} 50-68(1995)

\bibitem{HMSSZ} H. Hindano, J. Metcalf, H. Smith, C. Sogge, Y. Zhou:
{\em On Abstract Strichartz Estimates and The Strauss Conjecture For Nontrapping Obstacles},
arXiv:0805.1673v2.

\bibitem{IK} M. Ikawa:
{\em Decay of solutions of the wave equation in the exterior of several convex bodies}, 
Ann. Inst. Fourier(Grenoble), {\bf 38}, 1998, pp. 113-146.

\bibitem{Ka} L. V. Kapitanski:
{\em Cauchy problem for a semilinear wave equation ii}
Jour. Soviet Math., 62: 2746-2776, 1992.

\bibitem{KSS}  M. Keel, H. F. Smith, and C. D. Sogge:
{\em Almost global existence for some semilinear wave equations},
J. Anal. Math. {\bf 87} (2002), 265-279.

\bibitem{KT} M. Keel and T. Tao,
{\em Endpoint Strichartz estimates},
Amer. J. Math. {\bf 120} (1998), 955--980.


\bibitem{LS}{H. Lindblad and C. D. Sogge}:
{\em On existence and scattering with minimal regularity for semilinear
wave equations},
J. Funct. Anal. {\bf 130} (1995), 357--426.

\bibitem{LS2}{H. Lindblad and C. D. Sogge}:
{\em Long-time existence for small amplitude semilinear wave equations},
Amer. J. Math. {\bf 118} (1996), 1047--1135.

\bibitem{MSS}  G. Mockenhaupt, A. Seeger, and C. D. Sogge:
{\em Local smoothing of Fourier integral operators and Carleson-Sj\"{o}lin estimates}
J. Amer. Math. Soc., 6: 65-130, 1993.


\bibitem{NZ} S. Nonnenmacher, M. Zworski:
{\em Quantum Decay Rates in Chaotic Scattering},
Acta Math. {\bf 203}(2009), 149-233.

\bibitem{Pe} H. Pecher 
{\em Nonlinear small data scattering for the wave and Klein-Gordan equations}
Math. Z. ,185:261-270, 1984.

\bibitem{SS} H. F. Smith and C. D. Sogge:
{\em Global Strichartz estimates for nontrapping perturbations of
the Laplacian},
Comm. Partial Differential Equations {\bf 25}, (2000), 2171--2183.



\bibitem{So} C. D. Sogge:
Lectures on nonlinear wave equations, 2nd edition,
International Press, Boston, MA, 2008.

\bibitem{Str} Strichartz, R:
{\em Restriction of Fourier transform to quadratic surfaces and decay of sulutions to the wave equation}
Duke Math J. {\bf 44}(3), 705-714(1977)

\bibitem{T} D. Tataru:
{\em Strichartz Estimates for Operatiors with Nonsmooth Coefficients and the Nonlinear Wave Equation},
American Journal of Mathematics {\bf 122}(2000), 349-376

\bibitem{WZ} J. Wunsch and M. Zworski:
{\em Resolvent Estimates for Normally Hyperbolic Trapped Sets},
Ann. Henri Poincar\'{e} 12 (2011), 1349-1385

\bibitem{XY} X. Yu:
{\em Generalized Strichartz Estimates On Perturbed Wave Equation And Applications On Strauss Conjecture},
arXiv: 0905.0038v2.


\end{thebibliography}
\end{document}